\documentstyle[11pt]{article}

\font\tensym=msbm10
\font\sevensym=msbm7
\font\fivesym=msbm5

\newfam\symfam
\textfont\symfam=\tensym
\scriptfont\symfam=\sevensym
\scriptscriptfont\symfam=\fivesym
\def\sym{\fam\symfam\tensym}

\def\N{{\sym N}}

\newtheorem{example}{Example}[section]

\newtheorem{proposition}[example]{Proposition}

\newtheorem{lemma}[example]{Lemma}

\def\boxit#1#2{%
\setbox1=\hbox{\kern#1{#2}\kern#1}%
\dimen1=\ht1 \advance\dimen1 by #1 \dimen2=\dp1 \advance\dimen2 by #1
\setbox1=\hbox{\vrule height\dimen1 depth\dimen2\box1\vrule}%
\setbox1=\vbox{\hrule\box1\hrule}%
\advance\dimen1 by .4pt \ht1=\dimen1
\advance\dimen2 by .4pt \dp1=\dimen2  \box1\relax}
\catcode`\@=12

\newdimen\squaresize \squaresize=14pt
\newdimen\thickness \thickness=0.5pt

\def\square#1{\hbox{\vrule width \thickness
   \vbox to \squaresize{\hrule height \thickness\vss
      \hbox to \squaresize{\hss#1\hss}
   \vss\hrule height\thickness}
\unskip\vrule width \thickness}
\kern-\thickness}

\def\vsquare#1{\vbox{\square{$#1$}}\kern-\thickness}

\def\bbox#1{
\!\!
\begin{picture}(15,12)(1,2)
\put(0,0){\framebox(12,12){\hbox{\rm #1}}}
\end{picture}
}

\def\tildaa#1{{#1~\tilde{}}}

\def\young#1{
\vbox{\smallskip\offinterlineskip
\halign{&\vsquare{##}\cr #1}}}

\newdimen\unit
\def\o{$\scriptscriptstyle{{\rm o}}$}
\def\put(#1,#2)#3{\raise#2\unit\rlap{\kern#1\unit #3}\ignorespaces}

\def\g{{\unit=1mm
\hbox {\put(2,1.2){'}
       \put(1,2)\o
       \put(2,2)\o
       \put(3,2)\o
       \put(1.5,1)\o
       \put(2.5,1)\o
       \put(2,0)\o\
}\kern 3.5 \unit}}

\def\gfill{\leaders\hbox to 1.2em{\hss\g\hss}\hfill}

\def\tab#1|#2|{ {}_{#2}^{#1} }   

\def\ttab#1|#2|#3|#4|{
\vbox{\hbox{$\tab#1|#2|$} \hbox{$\tab#3|#4|$} } }

\def\Proof{\noindent \it Proof -- \rm}
\def\cqfd{\hspace{3.5mm} \hfill \vbox{\hrule height 3pt depth 2 pt width 2mm}
\bigskip}

\def\binomial#1#2{\left(\,\matrix{#1 \cr #2}\,\right)}

\def\vtr#1{\vrule height 0mm depth #1mm width 0mm}

\title{Performance evaluation of modulation methods: \\
       a combinatorial approach}

\author{D. Krob~\thanks{
\ LIAFA (CNRS) - Universit\'e Paris 7 - 2, place Jussieu - 75251 Paris
Cedex 05 -- France -- {\tt e-mail:} {\tt dk@liafa.jussieu.fr}
}\ ,\ \
E.A. Vassilieva~\thanks{
\ LIAFA (CNRS) - Universit\'e Paris 7 - 2, place Jussieu - 75251 Paris
Cedex 05 -- France -- {\tt e-mail:} {\tt katya@liafa.jussieu.fr}
}
}

\date{}

\begin{document}

\maketitle

\section{Introduction} \label{sec:intro}

Modulating a numeric signal corresponds to the fact of transforming
the digital sequence that represents it, into a wave form. Modulation
is therefore clearly a technique of main interest in a number of
ingeneering domains such as computer networks, mobile communications,
satellite transmissions, television diffusion, etc.

\smallskip
Due to their practical importance, modulation methods were therefore
widely studied in signal processing. The classical Proakis textbook
devotes for instance a full chapter to this subject (cf Chapter 5
of \cite{Pro}). One of the most important problem in this area is to
be able to design and to evaluate the performance characteristics of
the optimum receivers associated with a given modulation method. The
performance analyses that occur in such a context, reduce in
particular to the computation of various probability errors (see
again \cite{Pro} for more details).

\smallskip
Among the different modulation protocols, a rather important (in
practice) class consists in methods where the modulation references
(i.e. the wave forms associated with all possible digital sequences
of a given length) are also modulated and hence submitted to the
transmission noise. In this kind of situation, the demodulating
decision needs by consequence to account two noisy informations
(the transmitted signal and the transmitted references). The
computation of the probability errors appearing in such contexts,
involves therefore very often to compute the following type of
probability:
\begin{equation}\label{eq:1}
P(U < V)
=
P\left(\ \sum_{j=1}^{N}\, u_j^* u_j \ < \ \sum_{j=1}^{N}v_j^* v_j\,\right)
\ ,
\end{equation}
where the $u_i$ and $v_i$'s stand for independent complex Gaussian random
variables with arbitrary variances respectively denoted
\begin{equation}\label{eq:chidelta}
E[ u_j^*u_j] = \chi_j \, , \quad E[v_j^*v_j]=\delta_j
\end{equation}
for every $j \in [1,N]$ (see Section \ref{sec:barett} for more details).

\smallskip
The problem of computing explicitely this last probability was hence
studied by a number of researchers coming from signal processing (cf
\cite{Bar,Im,Pro,Tur}). The most interesting result in this direction
was obtained by Barett (cf \cite{Bar}) who proved that one can express
the probability defined by (\ref{eq:1}) as follows:
\begin{equation}\label{eq:2}
P(U < V)
=
\ \sum_{k=1}^{N} \,
\left( \
\prod_{j\not=k}   \, {1\over 1-\delta_k^{-1}\delta_j}\
\prod_{j=1}^{N} \, {1\over 1+\delta_k^{-1}\chi_j}
\ \right)
\ .
\end{equation}

\smallskip
It appears that this last formula can in fact be interpreted in a
purely combinatorial way, using Schur functions and Young tableaux
(see Section \ref{sec:barett} for the details). This new approach
already lead to the obtention of the first, both algorithmically
efficient and numerically stable, practical method for computing
the probability $P(U < V)$ (see again Section \ref{sec:barett}
or \cite{DKT,DKTV} for more informations).

\smallskip
In this paper, we continue the combinatorial study of Barett's
formula by showing that it is in fact highly connected with a
slight modification of a very classical bijection of Knuth (cf
\cite{Knu} or \cite{Ful} for a more recent presentation) between
pairs of Young tableaux of conjugated shapes and $\{0,1\}$-matrices.
These new considerations give us clearly a better understanding
of Barett's result. They also allowed us to obtain the first results
with respect to specializations of Barett's formula that were still
not known for the moment (see Section \ref{sec:specia}).


\section{Background}

\subsection{Partitions} \label{sec:part}

A {\it partition} is a finite nondecreasing sequence
$
\lambda = (\lambda_1, \lambda_2, \ldots, \lambda_m)
$
of positive integers. Graphically each such partition can be represented
by a diagram of $\lambda_1 + \ldots + \lambda_m$ boxes, called its
{\it Ferrers diagram}, whose $(m-i+1)$-th row contains $\lambda_i$ boxes
for every $1 \leq i \leq n$. The partition $\lambda = (2,2,4)$ is for
instance represented by the Ferrers diagram
$$
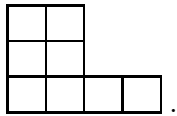
\young{ & \cr & \cr & & & \cr} \ .
$$
A given partition is then called the {\it shape} of the associated
Ferrers diagram.

\smallskip
Using this graphic representation, one can easily define the notion
of {\it conjugated partition}. The conjugated partition $\tildaa{\lambda}$
of a given partition $\lambda$ is indeed just the partition obtained
by reading the heights of the columns of the Ferrers diagram associated
with $\lambda$. One has here for instance $\tildaa{\lambda}=(1,1,3,3)$
when $\lambda = (2,2,4)$ as it can be seen on the previous picture.


\smallskip
When $\lambda$ is a partition whose Ferrers diagram is contained into
the square $(N^{N})$ with $N$ rows (of length $N$), one can also associate
with it its {\it complementary} partition, denoted by $\overline \lambda$,
which is the conjugate of the partition $\nu$ whose Ferrers diagram is the
complement (read from top to bottom) of the Ferrers diagram of $\lambda$
in the square $(N^{N})$. For instance, for $N = 6$ and $\lambda=(1,1,2,3)$,
we have $\nu=(3,4,5,5,6,6)$ and $\overline \lambda=(2,4,5,6,6,6)$ as it
can be checked on the following Figure~\ref{fig:1}. The Ferrers diagram
associated with $\lambda$ is here represented by the boxes filled with
$\bullet$ and the boxes filled with $\diamond$ correspond in the same
way to the partition $\nu$ (that can be obtained by computing the number
of such boxes per row) or to the complementary partition $\overline \lambda$
(that can be obtained by computing the number of such boxes per column).

\medskip

\begin{figure}[ht]
$$
\young{
\diamond & \diamond & \diamond & \diamond & \diamond & \diamond \cr
\diamond & \diamond & \diamond & \diamond & \diamond & \diamond \cr
\bullet  & \diamond & \diamond & \diamond & \diamond & \diamond \cr
\bullet  & \diamond & \diamond & \diamond & \diamond & \diamond \cr
\bullet  & \bullet  & \diamond & \diamond & \diamond & \diamond \cr
\bullet  & \bullet  & \bullet  & \diamond & \diamond & \diamond \cr
}
$$

\caption{\label{fig:1} Two complementary Young tableaux.}
\end{figure}

\smallskip
We will call {\it tabloid} of shape $\lambda$ the filling of a Ferrers
diagram of shape $\lambda$ with arbitrary positive integers. A filling
of the boxes of a Ferrers diagram of shape $\lambda$ with positive integers
is called a {\it Young tableau} (of shape $\lambda$) whenever the numbers
are weakly increasing along all rows and strictly increasing along all
columns. For example, the diagram
$$
\young{3 & 5\cr 2 & 2 \cr 1 & 1 & 1 & 4 \cr}
$$
is a Young tableau of shape $(2,2,4)$.

\smallskip
Let $X = \{\, x_i,\, 1 \leq i \leq n \, \}$ be a set of $n$ variables.
One associates then to any Young tableau $T$ filled by integers not
greater than $n$, a monomial $X^T$ defined as the product of the factors
$x_i$ for each entry $i$ of $T$. For the Young tableau $T$ of the above
example, one has for instance
$$
X^T = x_1^3\, x_2^2 \, x_3 \, x_4 \,x_5
$$
if we set $X = \{\, x_1, x_2, x_3, x_4, x_5 \}$. The {\it Schur function}
$s_\lambda(X)$ associated with the partition $\lambda$ is then defined as
the sum of the monomials $X^T$, for $T$ running over all Young tableaux of
shape $\lambda$, filled with numbers not greater than $n$. We recall that
each Schur function is a symmetric polynomial over $X$ and that the Schur
functions are a linear basis of the algebra of symmetric polynomials over
$X$ (for more informations on these questions, the reader should refer to
the classical textbook \cite{Mcd}).


\subsection{Gaussian polynomials and the $q-$Newton formula}
\label{sec:gauss}

Let $q$ be a variable. Then the expression
$$
\left[ n \right]_q
=
1+q+q^2+\ldots+q^{n-1}
=
\frac{1-q^n}{1-q}.
$$
is called the $q-$integer of order $n$ (this notation comes from the
fact that the specialization of a $q$-integer at $q = 1$ gives of course
the usual corresponding integer). We recall that the $q-$factorial is
then defined by
$$
\left[n!\right]_q
=
\left[1\right]_q \left[2\right]_q \ldots \left[n\right]_q
=
\frac{\displaystyle\prod_{i=1}^n (1-q^i)}{(1-q)^n}
\ .
$$
Finally the expression
$$
\binomial{n}{m}_q
=
\frac{\left[n!\right]_q}
{\left[m!\right]_q
\left[(n-m)!\right]_q}
$$
is known as the Gaussian polynomial of order $(n,m)$. It is clearly
the $q$-analogue of the usual binomial coefficient of the same order.
We refer to \cite{And} for more informations about Gaussian polynomials.

\smallskip
We will however recall the $q-$Newton formula (see \cite{And}):
\begin{equation}
\sum_{j=0}^N \
\binomial{N}{j}_q \,
(-1)^j\, q^{{j(j-1)}\over {2}}\, z^j
=
\ \prod_{k=0}^{N-1}\, (1-z\,q^{k}) \ .
\end{equation}
Note that the $q$-Newton formula specializes for $z=q$ to
\begin{equation} \label{equ:qnewton}
\sum_{j=0}^N \
\binomial{N}{j}_q\,
(-1)^j\, q^{{j(j+1)}\over {2}}
=
\ \prod_{k=1}^{N}\, (1-q^{k})\ .
\end{equation}


\subsection{Column bumping process and Knuth bijection}
\label{sec:knu}

The column bumping (or column-insertion) and row bumping (or row-insertion)
processes are classical dual constructions that provide algorithms to
transform a given word on the alphabet $\N$ of all integers into a Young
tableau. Here we will give some brief informations on the mechanism of
column bumping and present a famous bijection of Knuth, between pairs
of Young tableaux and $\{0,1\}$-matrices, (cf \cite{Knu}) which is,
being a variation on the well known Robinson-Schensted correspondence,
based on column bumping process in its construction. For a better
overview of the subject, we refer to \cite{Ful}.

\smallskip
The column bumping process is organized as follows. Take a positive
integer $x$ and a Young tableau $T$. Put $x$ in a new box at the top
of the first column if it is strictly larger than all the entries of
the column. If it is not the case, bump the lowest (i.e., the smallest)
entry in the column that is greater than or equal to $x$ and replace
it by $x$. Move the bumped entry to the top of the next column if
possible, or recursively bump one of the elements to the next column
otherwise. The process continues until the bumped entry can go at
the top of the next column, or until it becomes the only entry of a
new column. Note, that the bumping here takes place in a zig-zag path
that moves to the right, never moving up, and the result is always
another tableau. If the location of the box that is added is known,
the process can be reversed.

\smallskip
We are now in position to present the one-to-one correspondence (due to
Knuth) between matrices $M$ whose entries are zeros and ones (or equivalently
two-rowed arrays without repeated pairs) and pairs $(P,Q)$ of Young tableaux
with conjugated shapes. The construction of this bijection can be reflected
in the following steps.
\begin{enumerate}
\item Associate first with $M$ the array
$$
A
=
\left( \
\begin{array}{cccc}
u_1 & u_2 & \ldots \ldots & u_r \\
v_1 & v_2 & \ldots \ldots & v_r \\
\end{array}
\ \right),
$$
that consists of all the indices (classified in the lexicographic order)
corresponding to the $1$-entries of $M$ (all the entries $u_i$ of the
first row are therefore in weakly increasing order and one has moreover
$v_{i-1} \leq v_i$ in the second row whenever one has $u_{i-1} = u_i$
in the first).

\item Perform column bumping with all the variables $v_i$ of the second
row of the array $A$ begining from the first variable $v_1$ and moving
one by one to the very last variable $v_r$. The result is a Young tableau
$P$.

\item The second Young tableau $Q$ is just an encoding of the
order in which the first Young tableau $P$ was constructed on
the previous step. We first place the first element $u_1$ in the
(conjugated of the) first box that appeared during the column bumping
process that was used to construct $P$. The second element $u_2$
is placed in the same way in the box which is conjugated to the second
box that appeared in this process, etc.
\end{enumerate}
Applying the reversed column bumping process to the tableau $P$
and removing in the same time the corresponding boxes of the tableau
$Q$, allows us to reconstruct the initial array $A$ by writing
down in the order of their appearing the bumped out entries.


\subsection{Columns and their complements} \label{sec:colcomp}

In this subsection, we will finally pay some attention to {\it columns},
i.e. to Young tableaux of shape $1^k = (1,\ldots,1)$. The number $k$ of
1's is here equal to the length of the column. We will be only interested
by columns of length less than some positive number $N$, filled with integers
belonging to the set $\{1,\ldots,N\}$.

\smallskip
Introduce now some new notations. Let
$I = \{\,i_1 < i_2 < \ldots < i_l\, \}$
be a strictly increasing subsequence of $\{{1,\ldots,N}\}$. Denote then
by $c\left(I\right)$ the column of length $l$ filled with all integer
of $I$, increasing from bottom to top. We will also use in the sequel
the notation $c\left(I,J\right)$ to denote the column of length $l+m$
filled with the elements of the sets $I = \{i_1,i_2,\ldots,i_l\}$ and
$J = \{j_1,j_2,\ldots,j_m\},$ if the sequence
$$
i_1<i_2<\cdots<i_l<j_1<j_2<\cdots< j_m
$$
is a strictly increasing subsequence of $\{1,\ldots,N\}.$

\smallskip
Let $K = \{k_1 < k_2 < \ldots < k_t\}$ be a strictly increasing
subsequence of positive integers. The column $c\left(K \right)$ is then
called the {\it complement} (within $\{1,\ldots,N\}$) of the column
$c\left(I\right)$ if we have $K = \{1,\ldots, N\} \setminus I$. In
the sequel, this column will be denoted
$$
c\left(K\right) = \overline{c\left(I\right)} \ .
$$

\smallskip
Let us again take two sets $I = \{i_1,i_2,\ldots,i_l \}$ and
$J = \{j_1,j_2,\ldots,j_m\}$. We will say that the column $c\left(I\right)$
is {\it less or equal} than the column $c\left(J\right)$ and write
$c\left(I\right)\preceq c\left(I\right)$ if one has $m\leq n$ and
$i_k \leq j_k$ for every $1 \leq k \leq m$. In other words, a column
$c\left( I\right)$ is less or equal to a column $c\left(J\right)$ if
and only if one obtains a Young tableau when putting the column
associated with $J$ at the right of the column associated with $I$.


\section{Performance analysis of modulation protocols}

\subsection{Barret's formula} \label{sec:barett}

The analysis of many practical digital transmission systems involves
the computation of the probability that a given Hermitian quadratic
form in complex normal variates is negative. Numerous such examples
can be found in Proakis's standard textbook (cf \cite{Pro}). This
kind of problem appears in particular in the context of performance
analysis of classical demodulation protocols acting on modulated signals
transfered on noisy Gaussian channels.

\smallskip
A first expression for the probability that a given Hermitian quadratic
form in complex normal variates is negative, was first derived by Turin
(see \cite{Tur}) and used later by Barrett (see \cite{Bar}) to unveil a
closed form expression for this probability as a rational function of
the eigenvalues of the corresponding covariance matrix. Barett showed
indeed that the general problem discussed above can be reduced to the
study of the probability $P(U < V)$ presented in the first section of
this paper (cf formula (\ref{eq:1}) of section \ref{sec:intro}).

\smallskip
Barett gave also the explicit formula (\ref{eq:2}) that was, up to
this year, the best known approach from computing the probability
$P(U < V)$ from a practical point of view. Alternate methods involving
either direct contour integration of the associated characteristic
function along a carefully selected path so as to maximise numerical
stability, or algebraic manipulations like in \cite{Pro} (Annex B)
or \cite{Im}, provide other approaches involving numerical quadrature
of trigonometric functions.

\smallskip
All these methods lead however to algorithms that are not numerically
stable due to the presence of artificial singularities (such as the
singularities $\delta_i = \delta_j$ of Barett's formula (\ref{eq:2})).
It is therefore important to notice that the first efficient and stable
method for computing $P(U < V)$ was very recently proposed by Dornstetter,
Krob and Thibon (cf \cite{DKT} or section \ref{sec:cf}), based initially
on symmetric functions techniques. We recall below their algorithm for
the sake of completeness (cf \cite{DKT,DKTV} for all details).

\begin{itemize}
\item {\tt Step 1.} Consider the two polynomials defined by setting
$$
X(z)      = \ \prod_{i=1}^N \ (1-\chi_i\, z) \ ,
\qquad
\Delta(z) = \ \prod_{i=1}^N \ (1+\delta_i\, z) \ .
$$

\item {\tt Step 2.} Compute the unique polynomial $\pi$ of degree
less or equal $N\!-\!1$ such that one has
$$
\pi(z)\,X(z)+ \mu(z)\, \Delta(z) = 1
$$
where $\mu$ stands for some other polynomial of degree less or
equal to $N\!-\!1$.

\item {\tt Step 3.} Evaluate $P(U<V) = \pi(0)\ .$
\end{itemize}

\smallskip
\noindent
The algorithmic efficiency and the numerical stability of this
result comes then just from the fact that the second step of the
above method can be made using the generalized Euclidean algorithm
which is a very classical method which has the two above mentionned
properties.


\subsection{The combinatorial version of Barett's formula} \label{sec:cf}

Barett's formula in fact can be rewritten as a rational
fraction, i.e.
\begin{equation} \label{equ:sp}
P(U<V)
=
\frac{F(\chi,\delta)}{
\displaystyle\prod_{1 \leq i,j \leq N}(\chi_{i}+\delta_{j})
}\ ,
\end{equation}
where $F(\chi,\delta)$ is a symmetric polynomial with respect to the
$\chi_{i}$ and to the $\delta_{j}$. Moreover can be proved (cf
\cite{DKTV}) that $F(\chi,\delta)$ can be expressed in terms of
Schur functions, i.e.
\begin{equation} \label{eq:schur}
F(\chi,\delta)
=
\ \sum_{\lambda \subseteq (N^{N-1})} \
s_{(\lambda,N)}(\{\delta_1,\ldots,\delta_N\})\,
s_{\lambda^{\vee}}(\{\chi_1,\ldots,\chi_N\})
\ ,
\end{equation}
where $\lambda^{\vee}$ denotes the complement of the partition $(\lambda,N)$
within the rectangle $N^{N}$.

\smallskip
Note now that each monomial that appear in the right hand side of equation
(\ref{eq:schur}) can be obtained by taking the product of all the elements
of a square tableau of shape $N\times N$ consisting in two Young tableaux
of complementary shapes (i.e. as given by Figure \ref{fig:1} of Section
\ref{sec:part}) that respect the two following constraints:
\begin{itemize}
\item {\bf Condition S1}: the first tableau is only filled by variables
that belong to the alphabet $\delta = \{\,\delta_1,\ldots,\delta_N\,\}$
and {\it the length of its first row is equal to $N$},

\item {\bf Condition S2}: the second tableau is only filled by variables
that belong to the alphabet $\chi = \{\,\chi_1,\ldots,\chi_N\,\}$.
\end{itemize}
A typical example of such a combinatorial structure is given in
Figure \ref{fig:2}. Note that the first tableau is written here
in the usual way. On the other hand, the second tableau is organized
a bit differently: its rows (resp. its columns) are placed from top
to bottom (resp. from right to left) in the space corresponding to
the complement of the first tableau within the square $N\times N$.

\smallskip

\begin{figure}[ht]
$$
\young{
\chi_6   & \chi_5   & \chi_4   & \chi_3   & \chi_2   & \chi_1 \cr
\delta_4 & \chi_6   & \chi_5   & \chi_4   & \chi_2   & \chi_1 \cr
\delta_4 & \delta_5 & \delta_6 & \chi_4   & \chi_3   & \chi_2 \cr
\delta_3 & \delta_3 & \delta_5 & \chi_5   & \chi_4   & \chi_3 \cr
\delta_2 & \delta_2 & \delta_3 & \delta_4 & \chi_4   & \chi_3 \cr
\delta_1 & \delta_1 & \delta_2 & \delta_2 & \delta_2 & \delta_3 \cr
}
$$
\caption{\label{fig:2} A typical example of complementary filling
of a square tableau.}
\end{figure}

\smallskip
We will now proceed exploring the polynomial $F(\chi,\delta)$, involved
in formula (\ref{equ:sp}) and given by formula (\ref{eq:schur}), by
calculating the number $\alpha_N$ of all square $N \times N$ tableaux
filled as in the typical example of Figure~\ref{fig:2}. Knowing this last
integer will give us the exact algorithmic complexity of the formula
(\ref{equ:sp}). One should indeed just notice that $\alpha_N$ is equal
to the number of distinct monomials involved in $F(\chi,\delta)$, from
which one can easily deduce that the complexity of the computation
of $F(\chi,\delta)$ is exactly equal to $N^2\, \alpha_N$.

\smallskip
It appears unfortunately that $\alpha_N n= 2^{N^2-1}$, as proved in
the next result, which implies that formula (\ref{equ:sp}) can not
be used in practice as soon as $N$ grows. The combinatorial formula
(\ref{equ:sp}) is however absolutely not useless (from a theoretical
point of view) since it can be reformulated equivalently in the terms
of the algorithm given at the end of Section \ref{sec:barett}, which
is both practically very efficient (its complexity is quadratic as
Barett's formula) and numerically stable as already stated (cf
\cite{DKT,DKTV} for all details).

\begin{proposition}
The number $\alpha_N$ of square tableaux of shape $N\times N$ filled by
two complementary Young tableaux satisfying to conditions {\bf S1} and
{\bf S2} is given by the formula:
$$
\alpha_N = 2^{N^2-1} \ .
$$
\end{proposition}

\Proof Let us first notice that the conjunction of formulas
(\ref{eq:2}) and (\ref{equ:sp}) shows that one has:
\begin{equation} \label{eq:bar-sch}
\frac{F(\chi,\delta)}{
\displaystyle\prod_{1 \leq i,j \leq N}(\chi_{i}+\delta_{j})
}
=
\ \sum_{k=1}^{N} \,
\left(\
\prod_{j\not=k} \ {1\over 1-\delta_k^{-1}\delta_j}\
\prod_{j=1}^{N} \, {1\over 1+\delta_k^{-1}\chi_j}
\ \right)
\ .
\end{equation}

\smallskip
Let us begin by replacing everywhere $\chi_i$ and $\delta_i$ by $t^i$
in this last formula. This simple trick will allow us to avoid the
singularities of Barett's formula, corresponding to the situation
when some of the $\delta_i$'s collapse to a common value. Note that
this replacement of variables transforms the symmetric function
$F(\chi,\delta)$ into a polynomial $P(t)$ that provides the desired
number $\alpha_N$ when $t$ equal $1$. Therefore it is sufficient to
calculate $P(1)$ to get the value of $\alpha_N$.

\smallskip
Note now that formula (\ref{eq:bar-sch}) gives us immediately
the following expression for $P(t)$:
\begin{equation}
P(t)
=
\, \prod_{1\leq i,j\leq N} \,
\left(
t^i+t^j
\right) \,
\left( \ \,
\sum_{k=1}^N \,
\left( \,
\prod_{1 \leq j \neq k \leq N} \, \frac{1}{1-t^{j-k}}\
\prod_{j=1}^N \frac{1}{1+t^{j-k}} \
\right)\,
\right)\ .
\end{equation}

\smallskip
It appears that one can prove that the identity
\begin{equation}
\sum_{k=1}^N \,
\left(\
\prod_{1 \leq j\neq k \leq N}\, \frac{1}{1-t^{j-k}}\
\prod_{j=1}^N \frac{1}{1+t^{j-k}}
\ \right)
=
{1\over 2}
\end{equation}
holds for every $t$ (see Lemma \ref{lid} below). Hence one gets
$$
P(t)
=
{1\over 2}\,
\left(\
\prod_{1\leq i,j\leq N}\, (t^i+t^j)
\ \right) \ ,
$$
from which one can immediately conclude that $\alpha_N = P(1) = 2^{N^2-1}$.
Hence the proof of our proposition now reduces to the proof of the following
lemma.

\begin{lemma} \label{lid}
For every $t$, one has:
\begin{equation} \label{equ:lid}
\sum_{k=1}^N \,
\left( \
\prod_{1 \leq j\neq k \leq N} \, \frac{1}{1-t^{j-k}} \
\prod_{j=1}^N \, \frac{1}{1+t^{j-k}} \
\ \right)
=
{1\over 2}
\ .
\end{equation}
\end{lemma}
\Proof We will perform a number of equivalent transformations of
identity (\ref{equ:lid}) in order to reduce it into a classical
identity, which will finish our proof.

\smallskip
Taking first into account that
$$
\prod_{j=1}^N \ \frac{1}{1+t^{j-k}}
=
{1\over 2}\,
\left(\ \prod_{1 \leq j\neq k \leq N}\,
\frac{1}{1+t^{j-k}}
\ \right) \ ,
$$
we can rewrite equation (\ref{equ:lid}) in the equivalent way:
\begin{equation} \label{eq:iidd}
\sum_{k=1}^N\,
\left(\ \prod_{1 \leq j\neq k \leq N} \, \frac{1}{1-t^{j-k}} \ \right)
\left(\ \prod_{1 \leq j\neq k \leq N} \, \frac{1}{1+t^{j-k}} \ \right)
=
\sum_{k=1}^N\
\prod_{1 \leq j\neq k \leq N} \, \frac{1}{1-t^{2(j-k)}}
=
1
\ .
\end{equation}
Let us further develop the left hand side of the last equation.
We then get
$$
\sum_{k=1}^N \ \prod_{1 \leq j\neq k \leq N} \, \frac{1}{1-t^{2(j-k)}}
\ = \
\sum_{k=1}^N \ \prod_{j=1}^{k-1} \,\frac{1}{1-t^{2(j-k)}}\
         \prod_{j=k+1}^N \,\frac{1}{1-t^{2(j-k)}}
\vtr{5}
$$
$$
\ = \
\sum_{k=1}^N\ \prod_{j=1}^{k-1} \frac{1}{1-t^{(-2j)}}\
         \prod_{j=1}^{N-k} \frac{1}{1-t^{2j}}
\ = \
\sum_{k=1}^N\ \prod_{j=1}^{k-1}\, \frac{(-1)\,t^{2j}}{1-t^{2j}}\
         \prod_{j=1}^{N-k}\ \frac{1}{1-t^{2j}}
\vtr{5}
$$
$$
\qquad = \
\sum_{k=1}^N\ \frac{(-1)^{k-1}\,t^{2(1+2+\ldots+(k-1))}}
         {\displaystyle\prod_{j=1}^{k-1}\ (1-t^{2j})\
          \displaystyle\prod_{j=1}^{N-k}\ (1-t^{2j})}
\ = \
\sum_{k=1}^N\ \frac{(-1)^{k-1}\,t^{k\,(k-1)}}
         {\displaystyle\prod_{j=1}^{k-1}\ (1-t^{2j})\
          \displaystyle\prod_{j=1}^{N-k}\ (1-t^{2j})} \ .
$$
Substituting $t^2=q$ in the previous identity allows us therefore to
rewrite identity (\ref{eq:iidd}) into the following alternate form:
\begin{equation} \label{eq:iiiddd}
\sum_{k=1}^N \
\frac{(-1)^{k-1}\, q^{\frac{k(k-1)}{2}}}
{\displaystyle\prod_{j=1}^{k-1}\ (1-q^{j})
\displaystyle\prod_{j=1}^{N-k}\ (1-q^{j})}
=
1
.
\end{equation}
Multiplying now both parts on this last identity by the polynomial
$\left[(N\!-\!1)!\right]_q\, (1\!-\!q)^{N\!-\!1}$ and using the
definition of the Gaussian polynomials, we can rewrite identity
(\ref{eq:iiiddd}) as
$$
\left[(N\!-\!1)!\right]_q\, (1\!-\!q)^{N\!-\!1}
=
\
\sum_{k=1}^N \
\frac{(-1)^{k-1}\, q^{\frac{k(k-1)}{2}}\,
\left[(N\!-\!1)!\right]_q}
{\left[(k\!-\!1)!\right]_q\, \left[(N\!-\!k)!\right]_q}
=
\ \sum_{k=1}^N \
\binomial{\!N\!-\!1\!}{\!k\!-\!1\!}_q\, (-1)^{k-1}\, q^{\frac{k(k-1)}{2}}
\ .
$$
This last formula can therefore be equivalently rewritten as
$$
\prod_{i=1}^{N-1}\ (1-q^i)
=
\ \sum_{j=0}^{N-1} \
\binomial{\!N\!-\!1\!}{j}_q \,(-1)^j\, q^{\frac{j(j+1)}{2}}
\ ,
$$
which is exactly the well known $q-$Newton formula (see Section
\ref{sec:gauss} or \cite{And}). Hence the initial identity is true
since it is just a transformation of this last classical identity.
\cqfd


\section{A bijective proof of the combinatorial formula}

The previous proof gave us the desired number $\alpha_N$ of monomials
involved in $F(\chi,\delta)$ in a purely analytic way. It however did not
provide any insight, neither in the structure of $F(\chi,\delta)$, nor
in the simplicity of our result since the fact that $\alpha_N= 2^{N^2-1}$
is indeed clearly remarcable.

\smallskip
We will devote now this section to the construction of a bijective proof
of this last result. It will appear in fact that this construction will
also help us in studying a number of specializations of Barett's formula.
Hence our bijective proof will be rather interesting both from a theoretical
and a practical point of view.


\subsection{A more general structure}

In order to prove that $\alpha_N=2^{N^2-1}$ in a bijective way, we will
introduce a slightly generalized version of the combinatorial structures
that were involved in the description of $F(\chi,\delta)$. These new
combinatorial structures will just consist in the set, that we will denote
by ${\cal T}_N$, of all $N\times N$ squares divided into two complementary
Young tableaux (without any constraint on them) respectively filled by
elements of the alphabets $\delta$ and $\chi$. The following picture
shows two typical examples of an element of ${\cal T}_6$.

\smallskip

\begin{figure}[ht]
$$
\young{
\chi_6   & \chi_5   & \chi_4   & \chi_3   & \chi_2   & \chi_1 \cr
\delta_4 & \chi_6   & \chi_5   & \chi_4   & \chi_2   & \chi_1 \cr
\delta_4 & \delta_5 & \delta_6 & \chi_5   & \chi_2   & \chi_1 \cr
\delta_3 & \delta_3 & \delta_4 & \chi_6   & \chi_2   & \chi_1 \cr
\delta_2 & \delta_2 & \delta_2 & \delta_2 & \chi_2   & \chi_1 \cr
\delta_1 & \delta_1 & \delta_1 & \delta_1 & \delta_1 & \delta_1 \cr
}
\qquad
\qquad
\young{
\delta_6 & \chi_5   & \chi_4   & \chi_3   & \chi_2   & \chi_1 \cr
\delta_4 & \chi_6   & \chi_5   & \chi_4   & \chi_2   & \chi_1 \cr
\delta_4 & \delta_5 & \delta_6 & \chi_4   & \chi_3   & \chi_2 \cr
\delta_3 & \delta_3 & \delta_5 & \chi_5   & \chi_4   & \chi_3 \cr
\delta_2 & \delta_2 & \delta_3 & \delta_4 & \chi_4   & \chi_3 \cr
\delta_1 & \delta_1 & \delta_2 & \delta_2 & \delta_2 & \chi_4 \cr
}
$$
\caption{\label{fig:3} Two typical elements of ${\cal T}_6$\, .}
\end{figure}

\smallskip
\noindent
Note that the first tableau is again written in the usual way. On the
other hand, the second tableau is organized again differently: its rows
(resp. its columns) are placed from top to bottom (resp. from right
to left) in the space corresponding to the complement of the first
tableau within the square $N\times N$.

\smallskip
We will prove bijectively in the sequel that the cardinality of
${\cal T}_N$ is equal to $2^{N^2}$. This will immediately imply that
$\alpha_N = 2^{N^2-1}$ due to the fact that the number of elements
of ${\cal T}_N$ whose first tableau has a first row of length $N$
is clearly equal to the number of elements of ${\cal T}_N$ whose
second tableau has a first row of length $N$ (which corresponds
to the case where the first tableau has a first row of length
strictly less than $N$).

\smallskip
To get this last result, we will construct a bijection -- presented
in the next subsection -- between ${\cal T}_N$ and the set
${\cal M}_{N\times N}(\{0,1\})$ of all $\{0,1\}$-matrices of
size $N\times N$.


\subsection{Construction of the bijection}

We will now present our bijection between ${\cal M}_{N\times N}(\{0,1\})$
and ${\cal T}_N.$ Our construction is based on a slight variation of
the well known Knuth's bijection, presented in Section~\ref{sec:knu}.
We will see in the sequel that it has some deep and not obvious symmetry
properties that will be fundamental for highlighting Barett's formula
in a totally new way.

\smallskip
Let therefore $M$ be a matrix of ${\cal M}_{N\times N}(\{0,1\})$.
We associate then with $M$ the word $w(M)$ over the alphabet
$\{1,2,\ldots,N\}\times \{1,2,\dots,N\}$ defined as follows.
\begin{enumerate}
\item Construct first the $2$-row array $A_N$ which is equal to
the sequence of the $N^2$ pairs $(i,j)$ of
$\{1,2,\ldots,N\}\times \{1,2,\dots,N\}$ taken
in the lexicographic order, i.e.
$$
A_N
=
\left( \
\begin{array}{cccccccccc}
1 & \ldots & 1 & 2 & \ldots & 2 & \ldots \ldots & N & \ldots & N \\
1 & \ldots & N & 1 & \ldots & N & \ldots \ldots & 1 & \ldots & N \\
\end{array}
\ \right)
\ .
$$

\item Select then in this array all the entries that correspond to
the $1$'s of of $M$. We obtain then a word $w(M)$ on the alphabet
$\{1,2,\ldots,N\}\times \{1,2,\dots,N\}$ by reading all these entries
from left to right.
\end{enumerate}

\begin{example} \label{ex:w}
Let us consider the matrix
$$
M
=
\left( \
\begin{array}{ccc}
0 & 0 & 1 \\
1 & 0 & 0 \\
0 & 1 & 1 \\
\end{array}
\ \right) \ .
$$
Then one has
$$
A_3
=
\left( \
\begin{array}{ccccccccc}
1 & 1 & \bbox{1} & \bbox{2} & 2 & 2 & 3 & \bbox{3} & \bbox{3} \\
1 & 2 & \bbox{3} & \bbox{1} & 2 & 3 & 1 & \bbox{2} & \bbox{3}
\end{array}
\ \right)
$$
where we squared the entries associated with the $1$'s of $M$.
Hence we get
$$
w(M)
=
\left(
\begin{array}{c}
1 \\ 3
\end{array}
\right)
\
\left(
\begin{array}{c}
2 \\ 1
\end{array}
\right)
\
\left(
\begin{array}{c}
3 \\ 2
\end{array}
\right)
\
\left(
\begin{array}{c}
3 \\ 3
\end{array}
\right)
\ .
$$
\end{example}

We apply now Knuth's bijection to $w(M)$ in order to get two Young
tableaux
$$
\left(T_1, T_2\right)
$$
of conjugated shapes $\lambda_1$ and $\tildaa{\lambda_1}$. We will now
associate with the tableau $T_2$ a new tableau $\overline{T_2}$ of shape
$\overline{\lambda_1}$ (the complementary partition of $\lambda_1$ within
$N\times N$) that is constructed as follows. Note first that one defines
a unique tabloid $\overline{T_2}$ of shape $\overline{\lambda_1}$ by
asking (for every $i \in [1,N]$) that the $i$-th column of $\overline{T_2}$
consists exactly of all the letters of $\{\,1,\ldots,N\,\}$ that do not
appear in the $(N\!-\!i\!+\!1)$-th column of $T_2.$ It appears that this
tabloid is in fact a Young tableau.

\begin{proposition} \label{You}
The tabloid $\overline {T_2}$ is a Young tableau.
\end{proposition}

\Proof The proof will be made in three steps. In all lemmas that are
involved in this proof, we will use the notations and definitions of
Section~\ref{sec:colcomp}. The proof of the two first lemmas will be
found in the final version of this paper.

\begin{lemma} \label{cIJI}
Let $c\left(I,J\right)$ and $c\left(I\right)$ be the two columns
such that $c\left(I,J\right) \preceq c\left(I\right)$. Then, for
their complements $\overline{c\left(I,J\right)}$ and
$\overline{c\left(I\right)}$ holds the following
unequality:
$$
\overline{c\left(I\right)} \preceq
\overline{c\left(I,J\right)}
\ .
$$
\end{lemma}

\begin{lemma} \label{cIJ}
Let $c\left(I\right)$ and $c\left(J\right)$ be two columns of the
same length such that $c\left(I\right) \preceq c\left(J\right)$.
Then, for their complements $\overline{c\left(I\right)}$ and
$\overline{c\left(J\right)}$, holds the following inequality:
$$
\overline{c\left(J\right)} \preceq
\overline{c\left(I\right)}
\ .
$$
\end{lemma}

Proposition~\ref{You} is now an immediate consequence of the next
(and last) lemma.

\begin{lemma} \label{cIJK}
Let $c(I,J)$ and $c(K)$ be two columns that satisfy the inequality
$c\left(I,J\right) \preceq c\left(K\right)$. Suppose also that the
two subsets $I$ and $K$ of $\{1,\ldots,N\}$ have the same number of
elements. Then, for the complements $\overline{c\left(I,J\right)}$
and $\overline{c\left(K\right)}$ of the two above columns, one has:
$$
\overline{c\left(K\right)} \preceq
\overline{c\left(I,J\right)}
\ .
$$
\end{lemma}

\Proof The statement of our lemma can be easily obtained by applying
Lemma~\ref{cIJI} and Lemma~\ref{cIJ} in order to get the inequalities:
$$
c\left(I,J\right)
\preceq
c\left(I \right)
\preceq
c\left(K\right)
\ .
$$
This ends therefore both the proof of our lemma and of Proposition~\ref{You}.
\cqfd

\begin{example} \label{ex:t}
Let us continue Example~\ref{ex:w}. Knuth's bijection applied to the word
$w(M)$, gives the pair of tableaux
$$
\left(T_1,T_2\right)=
\left(\
\vcenter{
\hbox{
\young{
3        \cr
2        \cr
1 & 3    \cr
}\ \ ,
\young{
2          \cr
1 & 3 & 3     \cr
}
}
}
\
\right)
$$
of conjugated shapes $\lambda_1 =(1,1,2)$ and $\tildaa{\lambda_1}=(1,3)$.
The shape $\overline{\lambda_1}=(2,3)$, complementary to the shape $\lambda_1$
of the tableau $T_1$ within the square 3$\times$3, provides then the shape
of the tableau $\overline {T_2}$. Filling in its entries by taking (in the
reverse order) the complements within $\{1,2,3\}$ of the columns of $T2$,
we obtain
$$
\overline {T_2}=
\young{
2 & 2      \cr
1 & 1 & 3  \cr
}
\quad
.
$$
\end{example}

\smallskip
The pair $\left(T_1, \overline{T_2}\right)$ is then a pair of
complementary Young tableaux within the square $N\times N$. To
get an element of ${\cal T}_N$, it suffices now to associate with
each entry $i$ of $T_1$ (resp. of $T_2$) the letter $\delta_i$
(resp. $\chi_i$) of the alphabet $\delta$ (resp. $\chi$). The
element of ${\cal T}_N$ associated in such a way with the initial
matrix $M$ of ${\cal M}_{N\times N}(\{0,1\})$, will be denoted
by $\Phi(M)$ in the sequel.

\smallskip
Since the mapping $T_2\rightarrow \overline{T_2}$ is one to one,
it is now clear that we constructed in such a way a bijection
$\Phi$ between ${\cal M}_N(\{0,1\})$ and ${\cal T}_N.$ The problem
is now to explore the properties of this bijection in order to be
able to get some interesting enumerative consequences.

\begin{example}
Let us finish the previous example \ref{ex:t} which was itself a
continuation of Example \ref{ex:w}. The element of ${\cal T}_3$
which is associated with the pair $(T_1,\overline{T_2})$ is given
below:
$$
\Phi(M)
\ \
=
\ \
\young{
\delta_3 & \chi_2   & \chi_3 \cr
\delta_2 & \chi_2   & \chi_1 \cr
\delta_1 & \delta_3 & \chi_1 \cr
}
$$
where $M$ stands for the matrix introduced in Example \ref{ex:w}.
\end{example}


\subsection{Symmetry properties of our bijection}

We will present here a very strong symmetry property of our bijection
$\Phi$. To this purpose, we give first another method for constructing
it, presented below.

\begin{enumerate}
\item Construct again the two row array $A_N$ which is the sequence
of the $N^2$ pairs $(i,j)$ of $\{1,\ldots,N\}\!\times\!\{1,\ldots,N\}$
taken in the lexicographic order, i.e.
$$
A_N
=
\left( \
\begin{array}{cccccccccc}
1 & \ldots & 1 & 2 & \ldots & 2 & \ldots \ldots & N & \ldots & N \\
1 & \ldots & N & 1 & \ldots & N & \ldots \ldots & 1 & \ldots & N \\
\end{array}
\ \right) \ .
$$
Select then in $A_N$ all the pairs corresponding to the $1$'s of $M$.
We obtain then a first word $w_1(M)$ by reading the second component
of the selected entries.

\item Construct then the two row array $B_N$ which is equal to the sequence
of the $N^2$ pairs $(i,j)$ of $\{1,\ldots,N\}\!\times\!\{1,\ldots,N\}$
taken in the antilexicographic order (that is to say the lexicographic
order with respect to the second entry), i.e.
$$
B_N
=
\left( \
\begin{array}{cccccccccc}
1 & \ldots & N & 1 & \ldots & N & \ldots \ldots & 1 & \ldots & N \\
1 & \ldots & 1 & 2 & \ldots & 2 & \ldots \ldots & N & \ldots & N \\
\end{array}
\ \right)
\ .
$$
Select in this array all pairs corresponding to the $0$'s of $M$. We
obtain then a second word $w_2(M)$ by reading the first component of
the selected entries.
\end{enumerate}

\smallskip
\noindent
One construct then two Young tableaux $\left(T_1 ', T_2 '\right)$
by applying the column bumping process to the two previous words
$w_1(M)$ and $w_2(M)$. It appears that these tableaux are exactly
the two Young tableaux obtained by the bijection $\Phi$, constructed
in the previous subsection, applied to the matrix $M$.

\begin{example}
This example continues again Example~\ref{ex:t}. Since the first step
of the both ways to construct bijection $\Phi$ is the same, we will get
here
$$
A_3
=
\left( \
\begin{array}{ccccccccc}
1 & 1 & \bbox{1} & \bbox{2} & 2 & 2 & 3 & \bbox{3} & \bbox{3} \\
1 & 2 & \bbox{3} & \bbox{1} & 2 & 3 & 1 & \bbox{2} & \bbox{3}
\end{array}
\ \right) \ .
$$
For the second array, we have in the same way:
$$
B_3
=
\left( \
\begin{array}{ccccccccc}
\bbox{1} & 2 & \bbox{3} & \bbox{1} & \bbox{2} & 3 & 1 & \bbox{2} & 3 \\
\bbox{1} & 1 & \bbox{1} & \bbox{2} & \bbox{2} & 2 & 3 & \bbox{3} & 3
\end{array}
\ \right)
\ .
$$
Hence we get
$$
w_1(M)
=
\left(
3,1,2,3
\right)
\ ,
\quad
w_2 (M)
=
\left(
1,3,1,2,2
\right)
\ .
$$
The column bumping process applied to $w_1(M)$ and $w_2(M)$ gives us
then immediately the two following Young tableaux:
$$
\left(T_1',T_2'\right)
=
\left(\
\vcenter{
\hbox{
\young{
3        \cr
2        \cr
1 & 3    \cr
}\ \ ,
\young{
2 & 2        \cr
1 & 1 & 3     \cr
}
}
}
\
\right)
=
\left(T_1, \overline{T_2}\right)
\ .
$$
\end{example}

\medskip
We are now in position to state the following proposition which
expresses the main symmetry property of our construction.

\begin{proposition}
For every matrix $M$ of ${\cal M}_{N\times N}(\{0,1\})$, one has :
$$
\Phi(M) = \left(T_1 ', T_2 '\right)
$$
\end{proposition}

\Proof We will not give here the proof of this important result
sicne it is rather technical. It is just worthwhile to note that
our proof is based on the explicitation of the strong relations
that exist between the Greene's invariants of the two words
$w_1(M)$ and $w_2(M)$.
\cqfd

%
%


\subsection{Some specializations of Barett's formula} \label{sec:specia}

As a consequence of the bijection, we can get some interesting
combinatorial identities for several special cases of Barett's
formula. It is for instance obvious to see that our bijection
leads immediately to the identity
$$
F(\chi,\delta) + F(\delta,\chi)
=
\
\sum_{k=0}^{N^2} \
\binomial{N^2}{k}\,\delta^k\,\chi^{N^2-k}
$$
in the situation where one substitutes in the symmetric function $F$
all $\delta_i$ and $\chi_i$ by a single value respectively equal to
$\delta$ and $\chi$.

\smallskip
In a more interesting level, it is also possible to use our bijection
in order to get an explicit combinatorial interpretation (which was
still an open problem) of the coefficients of the polynomial of two
variables resulting from the substitution of the $k$ first variables
$\delta_i$ and $\chi_i$ by single values and of all last $N\!-\!k$
variables $\delta_j$ and $\chi_j$ by $1$. This interpretation is
however rather long to explain: it will therefore be only presented
in the final version of this paper.


\end{document}